\documentclass[12pt,twoside]{article}

\usepackage{mathtools}
\usepackage[colorlinks=true,linkcolor=blue,citecolor=blue,urlcolor=blue]{hyperref}
\usepackage{amsmath}
\usepackage{amsthm}
\usepackage[noabbrev]{cleveref}
\usepackage{amsfonts}
\usepackage{amssymb}
\usepackage{url}

\usepackage{graphicx}
\usepackage{float}
\usepackage[makeroom]{cancel}
\usepackage{framed}
\usepackage{scalerel}
\usepackage{color}
\usepackage{yfonts}

\usepackage{caption}

\def\dd{\hskip 1pt{\mathrm{d}}}

\def\bC{\mathbb{C}}
\def\bN{\mathbb{N}}

\def\diag{\mathop{\mathrm{diag}}}
\def\re{\mathop{\mathrm{Re}}}

\numberwithin{equation}{section}
\theoremstyle{plain}
\newtheorem{theorem}{Theorem}[section]

\newtheorem{proposition}[theorem]{Proposition}

\theoremstyle{definition}

\begin{document}

\title{\bf \large
A Reflection Formula for the Gaussian Hypergeometric Function of Matrix Argument
}

\author{{Donald Richards}{\hskip1pt}\thanks{Department of Statistics, Pennsylvania State University, University Park, PA 16802, U.S.A. E-mail address: \href{mailto:richards@stat.psu.edu}{richards@stat.psu.edu}.
\endgraf
\ $^\dag$Department of Mathematics \& Statistics, The College of New Jersey, Ewing, NJ 08628, U.S.A. E-mail address: \href{mailto:zheng@tcnj.edu}{zheng@tcnj.edu}.
\endgraf
\endgraf
\endgraf
\ {\it Running head}: Gaussian Hypergeometric Functions of Matrix Argument.} 
\ and {Qifu Zheng}{\hskip1pt}$^\dag$
}

\maketitle

\date

\begin{abstract}
We obtain a reflection formula for the Gaussian hypergeometric function of real symmetric matrix argument.  We also show that this result extends to the Gaussian hypergeometric function defined over the symmetric cones, and even to generalizations of the Gaussian hypergeometric function defined in terms of series of Jack polynomials.  Finally, we obtain a quadratic transformation formula for the Gaussian hypergeometric function of Hermitian $2 \times 2$ matrix argument.  
\end{abstract}


\normalsize

\parindent=25pt

\section{Introduction}
\label{introduction}
\setcounter{equation}{0}

The Gaussian hypergeometric series, typically denoted by ${}_2F_1(a,b;c;x)$, is well-known to satisfy a large number of linear (contiguous), quadratic, and nonlinear transformation formulas.  We refer to Andrews, {\it et al.} \cite[Chapter 3]{AAR}, Erd\'elyi, {\it et al.} \cite[Chapter II]{Erdelyi1}, or Olde Daalhuis \cite[Chapter 15]{Daalhuis} for extensive treatments of these formulas.  

In this paper, we are concerned with generalizations of a reflection formula, 
\begin{equation}
\label{reflection2f1}
{}_2F_1(a,b;a+b+1-c;1-x) = \frac{\Gamma(c-a)\Gamma(c-b)}{\Gamma(c)\Gamma(c-a-b)} \ {}_2F_1(a,b;c;x),
\end{equation}
valid for $0 < x < 1$, $-a \in \bN$, and $c, b-c+a+1 \notin \{0,-1,-2,\ldots,a+1\}$.  We extend the identity \eqref{reflection2f1} to the Gaussian hypergeometric function of matrix argument \cite{GrossRichards87,Herz,James,Muirhead}, to the Gaussian hypergeometric function defined on symmetric cones, and even to the Gaussian hypergeometric function defined in terms of Jack polynomials \cite{Faraut,Macdonald2}.  As a consequence, we establish a result that is stated without proof in the Digital Library of Mathematical Functions \cite[Eq.~(35.7.8)]{Richards}.  

We also consider the problem of deriving quadratic transformations for the Gaussian hypergeometric function of matrix argument, and we derive such a transformation for the hypergeometric function defined on the space of Hermitian $2 \times 2$ matrices.

\section{Remarks on the classical reflection formula \texorpdfstring{(\ref{reflection2f1})}{1.1}}
\label{classicalreflection}
\setcounter{equation}{0}

The reflection formula (\ref{reflection2f1}) is a special case of a well-known linear transformation of the Gaussian hypergeometric function; see \cite[p.~78, Theorem 2.3.2]{AAR} or \cite[Eq.~(15.8.7)]{Daalhuis}.  The formula is also stated indirectly by Whittaker and Watson \cite[p.~296, Miscellaneous Example 2]{WW}, who posed the problem of proving that if $-a \in \bN$ and $b, c \notin \bN$ then the ratio 
\begin{equation}
\label{Gaussian_hgf_ratio}
\frac{{}_2F_1(a,b;a+b+1-c;1-x)}{{}_2F_1(a,b;c;x)},
\end{equation}
is independent of $x$; further, the reader is asked to calculate the corresponding constant.  Once it has been proved that this ratio is independent of $x$, the constant of proportionality can be found by setting $x = 1$ and applying Gauss' formula, 
\begin{equation}
\label{Gauss}
{}_2F_1(a,b;c;1) = \frac{\Gamma(c) \, \Gamma(c-a-b)}{\Gamma(c-a) \, \Gamma(c-b)},
\end{equation}
for $\re(c-a-b) > 0$.  We remark that we were unable to infer, despite a close reading of \cite[Chapter XIV]{WW}, whether Whittaker and Watson had expected readers to solve this problem using a method alternative to the previously-mentioned linear transformation, and this motivated us to investigate several approaches to proving \eqref{reflection2f1}.  

For $a \in (0,1)$, $b = 1-a$, and $c = 1$, Berndt and Chan \cite{Berndt} showed that certain elliptic modular functions are expressible in terms of the inverse functions of ratios of the type \eqref{Gaussian_hgf_ratio}.  
For $n \in \bN$, Hannah \cite[p. 88]{Hannah} stated (\ref{reflection2f1}) in the form, 
$$
{}_2F_1(-n,b;c;1-x) = \frac{(c-b)_n}{(c)_n} \ {}_2F_1(-n,b;-n+b+1-c;x),
$$
$n \in \bN$, where $c, 1+b-c-n \notin \{0,-1,-2,\ldots,-n+1\}$, and remarked that a proof can be obtained by induction on $n$.  

Another direct but lengthy verification of \eqref{reflection2f1} can be obtained by expanding the left-hand side as a sum of powers, $(1-x)^j$, where $0 \le j \le a$; expanding each such power using the binomial theorem; reversing the order of summation; evaluating the inner sum using Gauss' summation theorem \eqref{Gauss}; and simplifying the resulting expression.  We were unable to locate in the extant literature an explicit presentation of this approach, but it is sufficiently straightforward that it is likely to have been folklore knowledge.  

We also remark that there is a certain finality to the reflection identity \eqref{reflection2f1} in the following sense.  Suppose that we are given the identity 
\begin{equation}
\label{reflection2f1_final}
{}_2F_1(\alpha,\beta;\gamma;1-x) = \tau \ {}_2F_1(a,b;c;x),
\end{equation}
$0 < x < 1$, for constants $\alpha$, $\beta$, $\gamma$, and $\tau$ that depend on $a$, $b$, and $c$, where $-a, -\alpha \in \bN$, and $\gamma$ and $c$ are such that both hypergeometric series are well-defined.  By repeatedly differentiating both sides of this identity with respect to $x$ and evaluating the outcome as $x \to 0$ and as $x \to 1$, i.e., by comparing coefficients of $x^k$ and also the coefficients of $(1-x)^k$, it can be deduced through relatively straightforward manipulations that $\alpha = a$, $\beta = b$, $\gamma = a+b+1-c$, and the constant $\tau$ can be evaluated in the usual way by letting $x \to 1$.

\section{A reflection formula for the Gaussian hypergeometric function of matrix argument}
\label{sec_Gaussianhgf}
\setcounter{equation}{0}

Another approach to establishing \eqref{reflection2f1} is by means of Euler's hypergeometric differential equation for the ${}_2F_1$ function.  As we now show, this approach generalizes to the Gaussian hypergeometric function of matrix argument; indeed, after an extensive investigation of potential approaches, it appears that the method of differential equations is the only approach that generalizes to higher-dimensional settings.  

First, we consider the case of $m \times m$ real symmetric matrix arguments, $X$.  For $\re(a) > \tfrac12(m-1)$, let 
$$
\Gamma_m(a) = \pi^{m(m-1)/4} \prod_{j=1}^m \Gamma\big(a - \tfrac12(j-1)\big)
$$
denote the \textit{multivariate gamma function}.  A \textit{partition} $\kappa = (k_1,\ldots,k_m)$ is a $m$-tuple of nonnegative integers $k_1,\ldots,k_m$ such that $k_1 \ge \cdots \ge k_m$; the \textit{weight} of $\kappa$ is $|\kappa| = k_1+\cdots+k_m$; the {\it partitional rising factorial} is 
$$
[a]_\kappa = \prod_{j=1}^m \big(a - \tfrac12(j-1)\big)_{k_j},
$$
where $(a)_k = a (a+1) \cdots (a+k-1)$ is the classical rising factorial; and we denote by $Z_\kappa(X)$ the corresponding {\it zonal polynomial} \cite{Faraut,GrossRichards87,James}.

For $a, b, c \in \bC$ and for $X$, a $m \times m$ real symmetric matrix, the {\it Gaussian hypergeometric function of matrix argument} is 
\begin{equation}
\label{Gaussianhgfoma}
{}_2F_1(a,b;c;X) = \sum_{k=0}^\infty \frac{1}{k!} \sum_{|\kappa| = k} 
\frac{[a]_\kappa \, [b]_\kappa}{[c]_\kappa} \ Z_\kappa(X),
\end{equation}
where the inner sum is over all partitions $\kappa$ of weight $k$.  The existence of this series requires that $[c]_\kappa \neq 0$ for all partitions $\kappa$, equivalently, $-c + 1 + \tfrac12(j-1) \notin \bN$ for all $j=1,\ldots,m$.  Denote by $x_1,\ldots,x_m$ the eigenvalues of $X$; if the series \eqref{Gaussianhgfoma} is non-terminating then it converges for all $m \times m$ real symmetric matrices $X$ such that $\|X\| := \max\{|x_1|,\ldots,|x_m|\} < 1$, and otherwise it converges for all such $X$ \cite{GrossRichards87}.  

The Gaussian hypergeometric function of matrix argument can also be defined by a generalized Euler-type integral over the cone of positive definite matrices \cite{GrossRichards87,Herz}.  Letting $I_m$ denote the $m \times m$ identity matrix then, as a consequence of those generalized Euler-type integrals, Herz \cite{Herz} generalized \eqref{Gauss} to: 
\begin{equation}
\label{Gauss_formula_hgfoma}
{}_2F_1(a,b;c;I_m) = \frac{\Gamma_m(c) \, \Gamma_m(c-a-b)}{\Gamma_m(c-a) \, \Gamma_m(c-b)},
\end{equation}
for $\re(c-a-b) > \tfrac12(m-1)$.  

The zonal polynomials $Z_\kappa(X)$ depend only on $x_1,\ldots,x_m$, the eigenvalues of $X$, and are symmetric functions of $x_1,\ldots,x_m$; therefore, the same holds for the hypergeometric function of matrix argument.  Muirhead \cite[p.~274, Theorem 7.5.5]{Muirhead} proved that the function ${}_2F_1(a,b;c;X)$ in \eqref{Gaussianhgfoma} is the unique solution of each of the partial differential equations 
\begin{multline}
\label{pdes_hgfoma}
x_i(1-x_i) \frac{\partial^2 F}{\partial x_i^2} + \Bigg[c - \tfrac12(m-1) - \big(a+b+1-\tfrac12(m-1)\big)x_i \\ 
+ \frac12 \sum_{\substack{j=1 \\ j \neq i}}^m \frac{x_i(1-x_i)}{x_i - x_j}\Bigg] \frac{\partial F}{\partial x_i} 
- \frac12 \sum_{\substack{j=1 \\ j \neq i}}^m \frac{x_j(1-x_j)}{x_i - x_j} \frac{\partial F}{\partial x_j} = abF,
\end{multline}
$i=1,\ldots,m$, subject to the conditions that:
\begin{itemize}
\setlength{\itemsep}{0.5pt}
\item[(i)] $F$ is a symmetric function of $x_1,\ldots,x_m$; 
\item[(ii)] $F$ is analytic at $X = 0$, i.e., $F$ is expressible in a neighborhood of $0$ as an infinite series of zonal polynomials $F(X) = \sum_{k=0}^\infty \sum_{|\kappa|=k} c_\kappa Z_\kappa(X)$, where the coefficients $c_\kappa$ do not depend on $m$; and 
\item[(iii)] $F(0) = 1$.  
\end{itemize}

\smallskip

The following result establishes a statement that was made without proof in \cite[Eq.~(35.7.9)]{Richards}.  

\smallskip

\begin{proposition}
\label{prop_reflection_hgfoma}
Suppose that $-a + 1 + \frac12(j-1) \in \bN$ for some $j=1,\ldots,m$.  Further, suppose that $-c + 1 + \tfrac12(j-1) \notin \bN$ and $-a - b + c - \tfrac12 (m - j) \notin \bN$ for all $j=1,\ldots,m$.  If $0 < X < I_m$, i.e., all eigenvalues of $X$ are in the interval $(0,1)$, then 
\begin{equation}
\label{reflection_hgfoma}
{}_2F_1\big(a,b;a+b+1-c + \tfrac12(m-1);I_m-X\big) = \frac{\Gamma_m(c-a)\Gamma_m(c-b)}{\Gamma_m(c)\Gamma_m(c-a-b)} \ {}_2F_1(a,b;c;X).
\end{equation}
\end{proposition}

\noindent{\it Proof}.  
Under the stated hypotheses on $b$ and $c$, the partitional rising factorials $[c]_\kappa$ and $[a+b+1-c + \tfrac12(m-1)]_\kappa$ are non-zero for all partitions $\kappa$; therefore, both sides of \eqref{reflection_hgfoma} are well-defined for all symmetric $m \times m$ matrices $X$ such that $\|X\| < 1$ and $\|I_m - X\| < 1$.  Moreover, under the assumption on $a$, both sides of \eqref{reflection_hgfoma} are terminating series, i.e., polynomials in $X$, hence are analytic at $X = 0$.  Moreover, by \eqref{Gauss_formula_hgfoma}, both sides are equal at $X = 0$.  

Replacing $x_i$ by $1-x_i$, $i=1,\ldots,m$, in each of the partial differential equations in the system \eqref{pdes_hgfoma}, we see that the outcome of this transformation is that $c$ is replaced by $a + b + 1 - c + \tfrac12(m-1)$.  Therefore the left- and right-hand sides of \eqref{reflection_hgfoma} each satisfy the system of differential equations in \eqref{pdes_hgfoma}.  Finally, by applying the uniqueness result of Muirhead \cite[Theorem 7.5.5]{Muirhead}, we find that the left- and right-hand sides of \eqref{reflection_hgfoma} coincide.
$\qed$

\bigskip

The reflection formula \eqref{reflection_hgfoma} extends to the hypergeometric function defined on the symmetric cones \cite{Faraut,GrossRichards87} or in terms of series of Jack polynomials \cite{Macdonald2}.  These generalizations satisfy a system of partial differential equations defined by means of the generalized Muirhead operators; {\it viz.}, for an arbitrary parameter $d > 0$, the system of partial differential equations is:
\begin{multline}
\label{pdes_hgfoma2}
x_i(1-x_i) \frac{\partial^2 F}{\partial x_i^2} + \Bigg[c - (m-1)d - \big(a+b+1- (m-1)d\big)x_i \\ 
+ d \sum_{\substack{j=1 \\ j \neq i}}^m \frac{x_i(1-x_i)}{x_i - x_j}\Bigg] \frac{\partial F}{\partial x_i} 
- d \sum_{\substack{j=1 \\ j \neq i}}^m \frac{x_j(1-x_j)}{x_i - x_j} \frac{\partial F}{\partial x_j} = abF,
\end{multline}
$i=1,\ldots,m$.  

Denote by $J_\kappa(x_1,\ldots,x_m;d)$ the corresponding Jack polynomials.  Then each of the equations \eqref{pdes_hgfoma2} has a common unique solution, denoted by ${}_2F_1(a,b;c;x_1,\ldots,x_n;d)$, subject to the conditions:  
\begin{itemize}
\setlength{\itemsep}{0.5pt}
\item[(iv)] $F$ is a symmetric function of $x_1,\ldots,x_m$; 
\item[(v)] $F$ is analytic at $(0,\ldots,0)$, i.e., $F$ is expressible in a neighborhood of the origin as a series of Jack polynomials $F(x_1,\ldots,x_m) = \sum_{k=0}^\infty \sum_{|\kappa|=k} c_\kappa J_\kappa(x_1,\ldots,x_m;d)$, where the coefficients $c_\kappa$ do not depend on $m$; and 
\item[(vi)] $F(0,\ldots,0) = 1$.  
\end{itemize}
For $\re(a) > (m-1)d$, let  
$$
\Gamma_m(a;d) = \pi^{m(m-1)d/2} \prod_{j=1}^m \Gamma\big(a - (j-1)d\big).
$$
denote the corresponding multivariate gamma function.  By applying the same argument used to establish Proposition \ref{prop_reflection_hgfoma}, we obtain the reflection formula, 
\begin{multline*}
{}_2F_1\big(a,b;a+b+1-c + (m-1)d;1-x_1,\ldots,1-x_m;d\big) \\
= \frac{\Gamma_m(c-a;d)\Gamma_m(c-b;d)}{\Gamma_m(c;d)\Gamma_m(c-a-b;d)} \ {}_2F_1(a,b;c;x_1,\ldots,x_m;d),
\end{multline*}
subject to the conditions that $-a + 1 + (j-1)d \in \bN$ for some $j=1,\ldots,m$, and $-c + 1 + (j-1)d \notin \bN$ and $-a - b + c - (m - j)d \notin \bN$ for all $j=1,\ldots,m$.

In closing this section, we remark that it would be interesting to extend to the Jack polynomial setting the approach described in the comments regarding \eqref{reflection2f1_final}.

\section{A quadratic transformation in the \texorpdfstring{$\boldsymbol{2 \times 2}$}{2by2} Hermitian case} 
\label{quadratictransformations}
\setcounter{equation}{0}

The problem of deriving nonlinear transformation formulas for the special functions of matrix argument was raised first by Herz \cite[p. 488]{Herz}, who noted the difficulty of deriving a quadratic transformation between the Bessel and confluent hypergeometric functions of matrix argument.  Such transformations are still generally unexplored and do not appear to follow from Euler integral representations, or manipulation of zonal polynomial or Jack polynomial series expansions.  It appears that such quadratic transformations will require analysis of the full monodromy group of the system of hypergeometric differential equations \eqref{pdes_hgfoma2}, and we note that Kor\'anyi \cite{Koranyi} studied a ``diagonal'' subgroup of that monodromy group and obtained a generalization of Kummer's twenty-four solutions for the classical Gaussian hypergeometric differential equation.  

A well-known quadratic transformation for the classical Gaussian hypergeometric function is 
\begin{equation}
\label{2F1quadratic}
{}_2F_1(\alpha,\alpha-\beta+\tfrac12;\beta+\tfrac12;t^2) = (1+t)^{-2\alpha} \ {}_2F_1\Big(\alpha,\beta;2\beta;\frac{4t}{(1+t)^2}\Big);
\end{equation}
see \cite[p. 176, Exercise 1(d)]{AAR} or \cite[Eq. (15.8.21)]{Daalhuis}.  We now illuminate the difficulties of deriving quadratic transformations in the matrix case by extending \eqref{2F1quadratic} to the hypergeometric function of \textit{Hermitian} matrix argument, i.e., defined on the space of $m \times m$ Hermitian matrices and corresponding to $d=1$ in \eqref{pdes_hgfoma2}.  

In the Hermitian case, the multivariate gamma function is
$$
\Gamma_m(a) = \pi^{m(m-1)/2} \prod_{j=1}^m \Gamma(a - j + 1),
$$
$\re(a) > m-1$; the partitional rising factorial is 
$$
[a]_\kappa = \prod_{j=1}^m \big(a-j+1)_{k_j};
$$
the zonal polynomial, $Z_\kappa$, is a multiple of the well-known Schur function $s_\kappa$ \cite{GrossRichards89,James}; and the Gaussian hypergeometric function of $X$, a $m \times m$ Hermitian matrix,  is 
\begin{equation}
\label{Hermitian2F1}
{}_2F_1(a,b;c;X) = \sum_{k=0}^\infty\frac{1}{k!} \sum_{|\kappa| = k} 
\frac{[a]_\kappa \, [b]_\kappa}{[c]_\kappa} Z_\kappa(X).
\end{equation}
The Gaussian hypergeometric function of $m \times m$ Hermitian matrix arguments, $X$ and $Y$, is defined as 
\begin{equation}
\label{Hermitian2F12arg}
{}_2F_1(a,b;c;X,Y) = \sum_{k=0}^\infty\frac{1}{k!} \sum_{|\kappa| = k} \frac{[a]_\kappa \, [b]_\kappa}{[c]_\kappa} \frac{Z_\kappa(X) \, Z_\kappa(Y)}{Z_\kappa(I_m)},
\end{equation}
In both series \eqref{Hermitian2F1} and \eqref{Hermitian2F12arg}, the parameter $c$ is such that $-c + j \notin \bN$ for all $j=1,\ldots,m$.  If the series \eqref{Hermitian2F1} is non-terminating then it converges for all $X$ such that $\|X\| < 1$; further, if \eqref{Hermitian2F12arg} is non-terminating then it converges if $\|X\| \cdot \|Y\| < 1$; see \cite[Theorem 6.3]{GrossRichards87}.  

Since the zonal polynomials depend only on the eigenvalues of their matrix arguments then we may assume, with no loss of generality, that $X = \diag(x_1,\ldots,x_m)$ and $Y = \diag(y_1,\ldots,y_m)$.  Define 
$$
V(X) = \prod_{1 \le i < j \le m} (x_i - x_j),
$$
and
$$
c_{2,1} = \beta_m^{-1} \frac{\prod_{i=1}^m (c-m+1)_{m-i}}{\prod_{i=1}^m (i-1)! \, (a-m+1)_{m-i} \, (b-m+1)_{m-i}}.
$$
The hypergeometric function of Hermitian matrix arguments can be expressed in terms of a determinant of classical hypergeometric functions.  Denote by $\det(a_{ij})$ the determinant of a $m \times m$ matrix with $(i,j)$th entry $a_{ij}$; then by \cite[Theorem 4.2]{GrossRichards89}, \cite{Khatri}, 
\begin{equation}
\label{Hermitian2F1reduction}
{}_2F_1(a,b;c;X,Y) = c_{2,1} \, \frac{\det\big({}_2F_1(a-m+1,b-m+1;c-m+1;x_iy_j)\big)}{V(X) V(Y)},
\end{equation}
where the ${}_2F_1$ functions on the right-hand side are the classical Gaussian hypergeometric functions, and L'Hospital's rule is to be applied if any of $x_1,\ldots,x_m$, or $y_1,\ldots,y_m$, coincide.  In particular, as in \cite[Eq. 5]{GrossRichards91}, by evaluating the limit as $Y \to I_m$, i.e., $y_1,\ldots,y_m \to 1$, it follows from \eqref{Hermitian2F12arg} and \eqref{Hermitian2F1reduction} that 
\begin{equation}
\label{Hermitian2F1reduction2}
{}_2F_1(a,b;c;X) = \frac{\det\big(x_i^{m-j} \, {}_2F_1(a-j+1,b-j+1;c-j+1;x_i)\big)}{V(X)}.
\end{equation}

Setting $(a,b;c) = (\alpha+m-1,\alpha-\beta+m-\tfrac12;\beta+m-\tfrac12)$, equivalently 
\begin{equation}
\label{abc_alphabeta}
(a-m+1,b-m+1;c-m+1) = (\alpha,\alpha-\beta+\tfrac12;\beta+\tfrac12),
\end{equation}
then we obtain 
\begin{align*}
c_{2,1}^{-1} \, V(X) V(Y) \ {}_2F_1(a,b;c;X,Y) &= \det\big({}_2F_1(a-m+1,b-m+1;c-m+1;x_iy_j)\big) \\
&= \det\big({}_2F_1(\alpha,\alpha-\beta+\tfrac12;\beta+\tfrac12;x_iy_j)\big).
\end{align*}
Applying \eqref{2F1quadratic}, we obtain 
\begin{multline}
\label{2F1X2Y2}
c_{2,1}^{-1} \, V(X^2) V(Y^2) \ {}_2F_1(a,b;c;X^2,Y^2) \\
= \det\left((1 + x_i y_j)^{-2\alpha} \ {}_2F_1\Big(\alpha,\beta;2\beta;\frac{4x_i y_j}{(1 + x_i y_j)^2}\Big)\right).
\end{multline}

From now on, we set $m = 2$.  Applying to the right-hand side of \eqref{2F1X2Y2} the $2 \times 2$ determinantal identity, 
\begin{equation}
\label{Hadamarddet}
\det(a_{ij} b_{ij}) = b_{11} b_{22} \det(a_{ij}) + a_{12} a_{21} \det(b_{ij}),
\end{equation}
we obtain 
\begin{align}
\label{2F1X2Y2expanded}
c_{2,1}^{-1} \, V(X^2) & V(Y^2) \ {}_2F_1(a,b;c;X^2,Y^2) \nonumber \\
= \ & \det\big((1 + x_i y_j)^{-2\alpha}\big) \cdot \prod_{j=1}^2 {}_2F_1\Big(\alpha,\beta;2\beta;\frac{4x_j y_j}{(1 + x_j y_j)^2}\Big) \nonumber \\
& + (1 + x_1 y_2)^{-2\alpha} (1 + x_2 y_1)^{-2\alpha} \cdot \det\left({}_2F_1\Big(\alpha,\beta;2\beta;\frac{4x_i y_j}{(1 + x_i y_j)^2}\Big)\right).
\end{align}

Next, we divide both sides of the latter equation by $V(Y^2) \equiv y_1^2 - y_2^2$, and let $Y \to I_2$, i.e., $y_1, y_2 \to 1$.  By applying L'Hospital's rule, we obtain 
$$
\lim_{Y \to I_2} \frac{\det\big((1 + x_i y_j)^{-2\alpha}\big)}{V(Y^2)} = \alpha (1 + x_1)^{-2\alpha - 1} (1 + x_2)^{-2\alpha - 1} (x_2 - x_1).
$$
Also, using the fact that 
$$
\frac{\partial}{\partial y_1} \dfrac{4x y_1}{(1 + x y_1)^2} \bigg|_{y_1 = 1} = \frac{4x(1-x)}{(1+x)^3},
$$
and the well-known formula, 
$$
\frac{\dd}{\dd y} \, {}_2F_1(a,b;c;y) = \frac{ab}{c} \, {}_2F_1(a+1,b+1;c+1;y),
$$
we obtain 
\begingroup
\addtolength{\jot}{1em}
\begin{align}
\label{temp_det1}
\lim_{Y \to I_2} & \frac{\det\left({}_2F_1\Big(\alpha,\beta;2\beta;\dfrac{4x_i y_j}{(1 + x_i y_j)^2}\Big)\right)}{V(Y^2)} \nonumber \\
&= \frac12 \lim_{y_1 \to 1} \frac{1}{y_1 - 1} 
\left|\begin{matrix}
{}_2F_1\Big(\alpha,\beta;2\beta;\dfrac{4x_1 y_1}{(1 + x_1 y_1)^2}\Big) & {}_2F_1\Big(\alpha,\beta;2\beta;\dfrac{4x_1}{(1 + x_1)^2}\Big) \\
& \\
{}_2F_1\Big(\alpha,\beta;2\beta;\dfrac{4x_2 y_1}{(1 + x_2 y_1)^2}\Big) & {}_2F_1\Big(\alpha,\beta;2\beta;\dfrac{4x_2}{(1 + x_2)^2}\Big)
\end{matrix}\right| \qquad\qquad \nonumber \\
& = \alpha \det(a_{ij} b_{ij}),
\end{align}
\endgroup
where, for $i,j=1,2$, 
\begin{align*}
a_{ij} &= \Big(\dfrac{1-x_i}{4(1+x_i)^2}\Big)^{2-j}, \\
b_{ij} &= \Big(\dfrac{4x_i}{(1+x_i)^2}\Big)^{2-j} \ {}_2F_1\Big(\alpha+2-j,\beta+2-j;2\beta+2-j;\dfrac{4x_i}{(1 + x_i)^2}\Big).
\end{align*}
Applying \eqref{Hadamarddet}, we find that the determinant $\det(a_{ij} b_{ij})$ in \eqref{temp_det1} equals 
\begin{align*}
& \dfrac{4x_1}{(1 + x_1)^2} \cdot \det\left({\hskip -2.2pt}\Big(\dfrac{1-x_i}{4(1+x_i)^2}\Big)^{2-j}\right)\cdot \prod_{j=1}^2 {}_2F_1\Big(\alpha+2-j,\beta+2-j;2\beta+2-j;\dfrac{4x_j}{(1 + x_j)^2}\Big) \\
& \, + \dfrac{1-x_2}{4(1+x_2)^2} \cdot \det\left(\Big(\dfrac{4x_i}{(1+x_i)^2}\Big)^{2-j} \, {}_2F_1\Big(\alpha+2-j,\beta+2-j;2\beta+2-j;\dfrac{4x_i}{(1 + x_i)^2}\Big)\right).
\end{align*}
Note that 
$$
\det\left(\Big(\dfrac{1-x_i}{4(1+x_i)^2}\Big)^{2-j}\right) = \frac{x_2 - x_1}{4(1+x_1)(1+x_2)}
$$
and, by \eqref{Hermitian2F1reduction2}, 
\begin{multline*}
\det\left(\Big(\dfrac{4x_i}{(1+x_i)^2}\Big)^{2-j} \ {}_2F_1\Big(\alpha+2-j,\beta+2-j;2\beta+2-j;\dfrac{4x_i}{(1 + x_i)^2}\Big)\right) \\
= V\big(4X(I_2+X)^{-2}\big)  \ {}_2F_1\big(\alpha+1,\beta+1;2\beta+1;4X(I_2+X)^{-2}\big).
\end{multline*}

The conclusion is that when both sides of \eqref{2F1X2Y2expanded} are divided by $V(Y^2)$ and then $Y \to I_2$, we obtain after some simplifications the result, 
\begin{align*}
c_{2,1}^{-1} \, V(X^2) & \ {}_2F_1(a,b;c;X^2) \\
= \ & \alpha (x_2 - x_1) \prod_{j=1}^2 (1 + x_j)^{-2\alpha - 1} \, {}_2F_1\Big(\alpha,\beta;2\beta;\frac{4x_j}{(1 + x_j)^2}\Big) \\
& + \alpha (1 + x_1)^{-2\alpha-2} (1 + x_2)^{-2\alpha-1} (x_2 - x_1) \\
& \quad \times \Bigg[\dfrac{x_1}{(1 + x_1)} \prod_{j=1}^2 {}_2F_1\Big(\alpha+2-j,\beta+2-j;2\beta+2-j;\dfrac{4x_j}{(1 + x_j)^2}\Big) \\
& \quad\qquad - \dfrac{(1-x_2)(1 - x_1 x_2)}{(1+x_2)^3} \ {}_2F_1\big(\alpha+1,\beta+1;2\beta+1;4X(I_2+X)^{-2}\big)\Bigg].
\end{align*}
Dividing both sides of the latter equation by $x_1 - x_2$ and using \eqref{abc_alphabeta} to substitute for $(a,b,c)$ in terms of $(\alpha,\beta)$, we obtain 
\begin{align*}
c_{2,1}^{-1} \, (x_1+x_2) & \ {}_2F_1(\alpha+1,\alpha-\beta+\tfrac32;\beta+\tfrac32;X^2) \\
= \ & - \alpha \prod_{j=1}^2 (1 + x_j)^{-2\alpha - 1} \, {}_2F_1\Big(\alpha,\beta;2\beta;\frac{4x_j}{(1 + x_j)^{2}}\Big) \\
& + \alpha (1 + x_1)^{-2\alpha-2} (1 + x_2)^{-2\alpha-1} \\
& \quad \times \Bigg[- \dfrac{x_1}{(1 + x_1)} \prod_{j=1}^2 {}_2F_1\Big(\alpha+2-j,\beta+2-j;2\beta+2-j;\frac{4x_j}{(1 + x_j)^{2}}\Big) \\
& \quad\qquad + \dfrac{(1-x_2)(1 - x_1 x_2)}{(1+x_2)^3} \ {}_2F_1\big(\alpha+1,\beta+1;2\beta+1;4X(I_2+X)^{-2}\big)\Bigg].
\end{align*}

\bigskip

\noindent
\textbf{Acknowledgments}.  The authors are grateful to the referee for helpful comments on the initial version of this article.

\vskip 0.3truein
\bibliographystyle{ims}

\end{document}